\documentclass[twoside,12pt, leqno]{article}
\usepackage{amsmath,amscd,amsthm,amssymb,amsxtra,latexsym,epsfig,epic,graphics}
\usepackage[matrix,arrow,curve]{xy}
\usepackage{graphicx}
\usepackage{diagrams}

\def\antiddot{\mathinner{\mkern1mu\raise1pt\vbox{\kern7pt\hbox{.}}\mkern2mu
        \raise4pt\hbox{.}\mkern2mu\raise7pt\hbox{.}\mkern1mu}}

\newcommand{\FF}{{\mathbb F}}

\newcommand{\PP}{{\mathbb P}}
\newcommand{\QQ}{{\mathbb Q}}

\newcommand{\ZZ}{{\mathbb Z}}

\newcommand{\s}{\mathcal}

\newcommand{\sF}{{\s F}}
\newcommand{\sG}{{\s G}}
\newcommand{\sH}{{\s H}}
\newcommand{\sI}{{\s I}}

\newcommand{\sL}{{\s L}}

\newcommand{\sO}{{\s O}}

\newcommand{\sZ}{{\s Z}}

\newcommand{\punkt}{\hspace{-.3ex}\raise.15ex\hbox to1ex{\Huge.}}

\DeclareMathOperator{\Pic}{Pic}

\DeclareMathOperator{\Spec}{Spec}

\DeclareMathOperator{\Hom}{Hom}
\DeclareMathOperator{\sHom}{\sH om}

\DeclareMathOperator{\pd}{pd}

\DeclareMathOperator{\rank}{rank}

\newtheorem{theorem}{Theorem}[section]

\newtheorem{conjecture}{Conjecture}[section]
\theoremstyle{definition}

\DeclareMathOperator{\Cliff}{Cliff}
\newcommand{\sExt}{{\mathcal E xt}}
\newcommand{\gM}{{\mathfrak M}}

\title{Computer aided Unirationality Proofs of Moduli Spaces}
\author{Frank-Olaf Schreyer}

\begin{document}
\maketitle
\begin{abstract}

We illustrate the use of the computer algebra system {\it Macaulay2} for  simplifications of classical unirationality proofs. We explicitly treat the moduli spaces of curves of genus g=10, 12 and 14.

\end{abstract}  

%
%

If a moduli space $M$ is unirational, then ultimately we would like to exhibit explicitly a dominating family,
which depends only on free parameters. Although this is possible along the lines of a known unirationality proof in principle, this is far beyond what computer algebra can do  today in most cases.

However, if we replace each generic free parameter with a random choice of an element in a ground field $\FF$, then in many cases the computation of a random element of the family defined
over $\FF$ is possible. 

This is particular interesting over a finite field $\FF$, because in that case there is no growth
of the coefficients in Gr\"obner basis computations. The  equidistributed probability  distribution on $\FF$ induces a probability measure on $M(\FF)$ via the algorithm, and opens up for us the possibility to investigate the moduli space experimentally.
With semi-continuity arguments we can sometimes deduce effectiveness of certain divisors on the moduli space. The computation with a single random example over a finite field can verify that the constructed family indeed dominates the desired moduli space. 
The key advantage in using Computer algebra instead of theoretical arguments is that random choices give smooth varieties almost certainly, while smoothness is always a difficult issue in any theoretical treatment.

In this note we illustrate this technique in a number of cases, and provide explicit  code for the computer algebra system {\it Macaulay2}. A more elaborate code is provided with the \textit{Macaulay2} package RandomCurves, which will be freely available from \cite{M2} and \cite{Sch}. 
I recommend running
the code  available from 
                   {www.math.uni-sb.de/ag/schreyer/home/forHandbook.m2} parallel with reading the code in the article. My favorite set-up of \textit{Macaulay2} is within emacs, which for example has syntax highlighting. That makes the code much more readable. Instructions how to use \textit{Macaulay2} with emacs
come with the installation package of \textit{Macaulay2}, see \cite{M2}. 

\noindent
{\bf Acknowledgements.}  I  thank the referees for their careful suggestions, and Florian Gei\ss \ and Hans-Christian Graf von Bothmer for valuable discussions.

\section{Random Plane Curves}

Let  $\gM_g$ denote the moduli space of curves of genus $g$.
For a general smooth projective curve $C$ of genus $g$  the Brill-Noether locus
$$W^r_d(C)= \{ L \in \Pic^d(C) \mid h^0(L) \ge r+1 \}$$
is non-empty  and smooth away from $W^{r+1}_d(C)$ of dimension $\rho$ if and only if 
the Brill-Noether number
$$\rho = \rho(g,d,r)=g-(r+1)(g-d+r) \ge 0,$$
equivalently, iff
$$d \ge ((r+1)(g-r)-g)/(r+1).$$
Moreover, $W^r_d(C)$ is connected in case $\rho>0$,
see \cite{ACGH}. The tangent space at the point $L \in W^r_d(C)\setminus W^{r+1}_d(C)$ is naturally  dual to the cokernel
of the Petri map
$$ H^0(C,L) \otimes H^0(C,\omega_C \otimes L^{-1})  \to H^0(C,\omega_C).$$
To prove the unirationality of $\gM_g$ we might prove that the Hilbert scheme of the corresponding models of $C$ in $\PP^r$ is unirational.

Severi's proof of the unirationality of $\gM_g$ for $g \le 10$ is based on the fact that a general curve up to genus 10
has a plane model with double points in general position, see \cite{Sev,AS}. An easy dimension count shows that this cannot be  true for $g\ge 11$. 

The algorithm that computes a random curve of genus $g= 10$ proceeds in four steps.
\begin{enumerate}
\item Compute the minimal degree $d$ such that $\rho(g,2,d)\ge 0$. 
\item Choose a scheme  $\Delta$ of $\delta={d-1 \choose 2}-g$ distinct points in $\PP^2$. 
\item Choose a random element in $f  \in H^0(\PP^2, \sI^2_\Delta(d))$.
\item Certify that $f$ defines an absolutely irreducible $\delta$-nodal curve.
\end{enumerate}
\noindent
In case $g=10$ we have $d=9$ and $\delta=18$.

The \textit{Macaulay2} code for these steps looks at follows. The code is, I believe, fairly self explaining, 
since the \textit{Macaulay2} language is close to standard mathematical notation. For further explanations I refer to the online documentation  
                   {http://www.math.uiuc.edu/Macaulay2/}.  \medskip

\noindent {\sc Step 1.} 
\begin{verbatim}
  i1 :  -- function that computes the minimal d: 
        dmin=(r,g)->ceiling(((r+1)*(g+r)-g)/(r+1)) 
  i2 :  g=10,r=2,d=dmin(r,g),delta=binomial(d-1,2)-g
  i3 :  rho=g-(r+1)*(g+r-d)
  o3 =  1  
\end{verbatim}

\noindent {\sc Step 2}. We specify $\delta$ points in $\PP^2$:
\begin{verbatim}
  i4 :  FF=ZZ/10007 -- a finite ground field
  i5 :  S=FF[x_0..x_2] -- the homogeneous coordinate ring of P2
  i6 :  -- pick a list of delta ideals of random points: 
        points=apply(delta,i->ideal(random(1,S),random(1,S)))
\end{verbatim}
Remark: The \textit{Macaulay2} function {\tt"random(n,S)"} picks a random element of degree $n$ in a graded ring $S$.
\begin{verbatim} 
  i7 :  IDelta=intersect points
  i8 :  betti res IDelta
              0 1 2
  o8 = total: 1 4 3
           0: 1 . .
           1: . . .
           2: . . .
           3: . . .
           4: . 3 .
           5: . 1 3      
  o8 = BettiTally
\end{verbatim}
For a very small finite field we cannot find $\delta$ distinct points defined over $\FF$. 
Even for $p=101$ the procedure will not  choose distinct points in $\PP^2(\FF_{101})$ in about $1.5 \%$ of all cases. For $p$ as large as $10007$ the failure probability is neglectable. 

A method to get points which works for small finite  fields is the following:
Since points have codimension 2 in $\PP^2$, we might use
the Hilbert-Burch matrix, see \cite{Eis}, Theorem 20.15. We recall the Hilbert-Burch theorem in its most useful special case:

\begin{theorem} \label{HB}
Let $R$ be a local (or graded) regular noetherian ring. The minimal free resolution of $R/I$ for a Cohen-Macaulay ideal $I \subset R$ of codimension 2 has shape
$$ 0 \leftarrow R/I  \leftarrow R \leftarrow R^{n+1} \leftarrow R^n \leftarrow 0$$
where the map $R \leftarrow R^{n+1}$ is given by the $n+1$ maximal minors of the matrix  $R^{n+1} \leftarrow R^n$.
Conversely, given a matrix $(a_{ij})_{i=0,\dots,n,j=1,\ldots,n}$, whose maximal minors have no common factor, then the ideal generated by the maximal minors is Cohen-Macaulay of codimension 2 and the corresponding complex is exact.
\end{theorem}
\noindent
Note  that 
$(m_0,\ldots,m_n)(a_{ij})=0$ for $m_k=(-1)^k \det (a_{ij})_{i\not=k}$
follows by expanding the determinants
$$0=\det 
\begin{pmatrix} 
a_{01} & \ldots & a_{0n} &a_{0i} \cr
\vdots & & \vdots & \vdots \cr
a_{n1} & \ldots & a_{nn} &a_{ni} \cr
\end{pmatrix}
$$
of matrices with a repeated column with respect to the last column.
The Hilbert-Burch Theorem is the basic reason why the Hilbert scheme of finite length subschemes
on a smooth surface is smooth, see \cite{Fog}.\medskip

\noindent
We continue with our \textit{Macaulay2}-program
\begin{verbatim}
  i9 :  M=random(S^{3:-4,1:-5},S^{3:-6}); -- a Hilbert-Burch matrix;
  i10 :  betti M
  i11 :  IDelta=minors(3,M); -- and its minors
  i12 :  betti res IDelta 
   
               0 1 2
  o12 = total: 1 4 3
            0: 1 . .
            1: . . .
            2: . . .
            3: . . .
            4: . 3 .
            5: . 1 3
  o12 : BettiTally
  i13 : degree IDelta==delta
  o13 = true
\end{verbatim}

\noindent {\sc Step 3.} Compute curves with double points in Delta:
\begin{verbatim}
  i14 : -- the saturation of the ideal IDelta^2 contains
        -- all equations which vanish double at Delta:
        J=saturate(IDelta^2); 
  i15 : betti J 
   
               0  1
  o15 = total: 1 10
            0: 1  .
            1: .  .
            2: .  .
            3: .  .
            4: .  .
            5: .  .
            6: .  .
            7: .  .
            8: .  1
            9: .  9
  o15 : BettiTally
\end{verbatim}
As expected, there is precisely one curve with multiplicity 2 at every point of $\Delta$.
\begin{verbatim}
  i16 : IC=ideal(gens J*random(source gens J, S^{-d}))
  i17 : degree IC
  o17 = 9
\end{verbatim}

\noindent {\sc Step 4.} 
To certify that we have indeed obtained a $\delta$-nodal curve $C$, it suffices to prove that the singular locus is reduced of degree $\delta$ because only plane curves with at most nodes have a reduced singular locus.
\begin{verbatim}
  i18 : singC=ideal jacobian IC + IC;
  i19 : codim singC == 2 and degree singC == delta
  o19 = true
  i20 : -- remove primary component at the irrelevant ideal:
        singCs=saturate(singC); 
  i21 :   betti singCs, betti singC 

                0 1         0 1
  o21 = (total: 1 4, total: 1 3)
             0: 1 .      0: 1 .
             1: . .      1: . .
             2: . .      2: . .
             3: . .      3: . .
             4: . 3      4: . .
             5: . 1      5: . .
                         6: . .
                         7: . 3

  o21 : Sequence
  i22 : codim (minors(2,jacobian singCs)+singCs) == 3
  o22 = true
\end{verbatim}
We verify that $C$ is geometrically irreducible from our information so far. Indeed, if the curve decomposes as 
$C_1 \cup C_2$ with $(\deg C_1, \deg C_2)=(a,b)$, say $a < b$, and $a+b=9$, then 
the intersection points $ C_1 \cap C_2$ form a subset of $\Delta$. So $(a,b)=(1,8)$ or $(2,7)$
is excluded because $I_\Delta$ is generated in degree $6$. The case $(4,5)$ is excluded since $20 > 18=\delta$
and $(3,6)$ is excluded because $\Delta$ is not a complete intersection. \medskip

\noindent
Finally, we conclude from this computation that $\gM_{10}$ is unirational over  $\QQ$.
We first note that the computations in the finite prime field $\FF$ may be viewed as the reduction mod p of the analogous computations for a curve defined over an open part of $\Spec \ZZ$.
By semi-continuity, the curve over $\QQ$ is $\delta$-nodal as well, hence has geometric genus 10, which proves that we have a rational map
$\mathbb A^n \dasharrow \gM_{10}$ defined over ${\QQ}$ for a suitable $n$. 
More precisely, we have a correspondence
$$
\begin{xy}
\xymatrix{ \operatorname{Hilb}_\delta(\PP^2) &
H=\{ ( C',\Delta) \mid C' \in \PP H^0(\PP^2, \sI^2_\Delta(d))^* \hbox{ is $\delta$-nodal } \} \ar[l] \ar[d]  \\\
&\gM_g}
\end{xy}
$$
For $g=10$ the left arrow gives a birational map of a component $H'$ of $H$ onto a dense open subset
of $\operatorname{Hilb}_\delta(\PP^2)$, because
$\chi( \sI^2_\Delta(9))={11 \choose 2} -3\delta=1$, and 
$h^0( \sI^2_\Delta(9))=1$ holds for our specific point $\Delta \in \operatorname{Hilb}_\delta(\PP^2)$.

The downward arrow factors through the universal $\mathcal W^2_d \subset \mathcal Pic^d_g \to \gM_{g}$, which has codimension at most $(2+1)(g+2-d)$ in the universal Picard variety over $\gM_g$ at every point.
A fiber of $H \to \mathcal W^2_d$ over a point $(C,L)$ with $h^0(C,L)=3$ is precisely the $PGL(3)$ orbit.
Thus to prove that $H' \to \gM_{10}$ dominates, it remains to certify that the fiber of
 $\mathcal W^2_d \to \gM_{10}$ over our specific point $C$ has expected dimension $\rho$ at
 the specific line bundle $L=\eta^* \sO_{\PP^2}(1) \in W^2_d(C)$ where $\eta:C \to C'\subset \PP^2$ denotes the normalization map.

By Riemann-Roch and adjunction $h^0(C,L)=r+1=3$ holds because $h^1(C,L)=h^0(C,\omega_C \otimes L^{-1})=h^0(\PP^2, \sI_\Delta(d-4))=3$. Moreover, the 
Petri map can be identified with
$$H^0(\PP^2,\sI_\Delta(d-4)) \otimes H^0(\PP^2, \sO(1)) \to H^0(\PP^2,\sI_\Delta(d-3))$$
which is an injection, since there is no linear relation among the three quintic generators of $I_\Delta$ by   the shape of the Hilbert-Burch matrix, compare also with the output line o8 and o12  of the \textit{Macaulay2} program above. Thus $W^2_d(C)$ is actually smooth of dimension $\rho$ in $L$ as desired. \medskip

It is easy to transform the code above into a function which chooses randomly $\delta$-nodal curves of degree $d$ provided the expected $h^0(\PP^2,\sI_\Delta(d))={d+2\choose 2}-3\delta >0$. An implementation can be found in the \textit{Macaulay2} package RandomCurves.

\section{Searching}

Over a finite field $\FF$ we might find points $C \in M(\FF)$  of a moduli space, provided $M$ is dominated by a variety $H$ of fairly low codimension in a unirational variety $G$. Indeed,
if $H$ is absolutely irreducible, then the proportion of  $\FF$-rational points is approximately
$$\frac{\mid H(\FF) \mid}{\mid G(\FF) \mid } \sim q^{-c}$$
where $q= \mid \FF \mid$ and $c$ is the codimension of $H \subset G$ by the Weil formulas.
If we can decide $p \notin H$ fast enough computationally, then we might be able to find points in $H(\FF_q)$ for small $q$ by picking points at random in $G(\FF_q)$ and testing $p \in H$. I will illustrate this technique by searching for plane models of random curves of genus 11.
Of course, to get just a random curve of genus 11, it is much better  to use Chang and Ran's unirational parameterization of $\gM_{11}$ via space curves, as indicated in the next section. 

This time we will use a bit more of the \textit{Macaulay2} syntax. In particular we will illustrate the use of method functions with options.

\begin{verbatim}
i23: randomDistinctPlanePoints = method(TypicalValue=> Ideal); 
     -- create the ideal of k random points 
    -- via a Hilbert-Burch matrix 
i24: randomDistinctPlanePoints (ZZ,PolynomialRing) := (k,S) -> (
     if dim S =!= 3 then error "no polynomial ring in 3 variables"; 
     -- numerical data for the Hilbert-Burch matrix
     n := ceiling((-3+sqrt(9.0+8*k))/2); 
     eps := k - binomial(n+1,2);
     a := n+1-eps; 
     b := n-2*eps; 
     distinct := false; 
     while not distinct do (
         -- the Hilbert-Burch matrix is B 
         B := if b >= 0 then random(S^a,S^{b:-1,eps:-2}) 
                 else random(S^{a:0,-b:-1}, S^{eps:-2});
     I := minors(rank source B, B); 
     distinct = distinctPlanePoints I);
     return I);
i25: distinctPlanePoints=method(TypicalValue=>Boolean);
i26: distinctPlanePoints(Ideal):= I-> (
           dim I==1 and dim (I+minors(2,jacobian I))<=0)
i27: distinctPlanePoints(List):= L->(
     -- for a List of ideals of points 
     -- check whether they have some point in common.
        degree intersect L == sum(L,I->degree I))
\end{verbatim}
The function {\tt randomDistinctPlanePoints} will return an ideal of a set of distinct points which we use
in our search below. 

The minimal degree of plane models of a general curve of genus $g=11$ is $d=10$, and we expect that a general model will have $\delta={d-1 \choose 2}-11=25$ ordinary double points. Since
$\chi(\PP^2, \sI^2_\Delta(10))=-9$ we expect that $\Delta$ with
$h^0(\PP^2, \sI^2_\Delta(10))=1$ from a codimension 10 subfamily of $Hilb_{25}(\PP^2)$. 
We can improve our odds if we look for models with 2 triple points.
Since triple points occur in codimension 1 and $\rho(11,2,10)=2$ we expect that the family of plane curves of degree 10 with two triple points $p_1,p_2$ and 19 double points will dominate $\gM_{11}$.
Since $\chi(\PP^2, \sI^2_\Delta \otimes \sI^3_{p_1}\otimes \sI^3_{p_2}(10))={10+2 \choose 2}-3\cdot 19-6\cdot 2=-3$,  we expect a
search for points in a subfamily of codimension 4. Thus the search function below  
will have
an average running time of order $O(q^4)$ with respect to the number of elements of $\FF=\FF_q$. So this will be only feasible for very small finite fields. 

\begin{verbatim}
i28: searchPlaneGenus11Curve=method(Options=>{Attempts=>infinity})
    -- search a plane curve of degree 10,geometric genus 11 with
    -- two triple points and 19 double points
i29: searchPlaneGenus11Curve PolynomialRing := opt -> S -> (
     I1 := ideal(S_0,S_1); 
     I2 := ideal(S_1,S_2); 
     IDelta := ideal 0_S; J := ideal 0_S;
     h := -3; attempt := 0; 
     while h <= 0 and attempt < opt.Attempts do (
           IDelta = randomDistinctPlanePoints(19,S); 
           while not distinctPlanePoints({I1,I2,IDelta}) do (
                IDelta = randomDistinctPlanePoints(19,S)); 
          J = intersect(I1^3,I2^3, saturate(IDelta^2));
          h = (tally (degrees gens truncate(10,J))_1)_{10}; 
          attempt = attempt +1);
     print attempt; 
     if attempt >= opt.Attempts then return null; 
     f := ideal 0_S; 
     gJ := gens J; 
     if h === 1 then f = ideal gJ_{0} 
     else while f == 0 do f=ideal(gJ*random(source gJ, S^{ -10})); 
     singf := ideal singularLocus f; 
     doublePoints := saturate(singf, I1*I2); 
     if degree doublePoints == 19 and (
         f31 := ideal contract(S_2^7, gens f);
        dim singularLocus f31 == 1) 
     and (
        f32 := ideal contract(S_0^7, gens f);
        dim singularLocus f32 == 1) 
    then return f else 
    searchPlaneGenus11Curve(S, Attempts => opt.Attempts-attempt));
\end{verbatim}
After these preparations the commands below will return a desired curve within a few minutes.
\begin{verbatim}
  i30 : p=5;FF=ZZ/p -- a finite ground field
  i32 : S=FF[x_0..x_2]
  i33 : setRandomSeed("alpha")
  i34 : C=searchPlaneGenus11Curve(S,Attempts=>2*p^4)
  18
  432
        -- used 48.132 seconds
  o34 : Ideal of S
\end{verbatim}

\section{Space Curves}\label{spaceCurves}

The proof of the unirationality of the moduli space $\gM_g$  for $g=11,12$ and $13$ by Sernesi and Chang-Ran is based on models of these curves in $\PP^3$, see \cite{Ser, CR}.

 Suppose $C \subset \PP^3$ is a Cohen-Macaulay curve with ideal sheaf $\sI_C$. The Hartshorne-Rao module
 $$ M = \sum_{n \in \ZZ} H^1(\PP^3, \sI_C(n)),$$
which has finite length and measures the deviation for $C$ to be projectively normal, plays an important role in liaison theory of curves in $\PP^3$. We briefly recall the basic facts.

Let $S=\FF[x_0,\ldots,x_3]$ and $S_C=S/I_C$ denote the homogeneous coordinate ring of $\PP^3$ and $C \subset \PP^3$ respectively.
By the Auslander-Buchsbaum-Serre formula, \cite{Eis} Theorem 19.9, $S_C$ has projective dimension $\pd_S S_C \le 3$ as an $S$-module. 
Thus the minimal free resolution has shape
$$ 0 \leftarrow S_C \leftarrow S \leftarrow F_1  \leftarrow F_2  \leftarrow F_3  \leftarrow 0 $$
with free graded modules $F_i=\oplus S(-j)^{\beta_{ij}}$.

The sheafified kernel
$\sG=ker(\widetilde F_1 \to \sO_{\PP^3})$ is always a vector bundle,
and
$$ 0\leftarrow\sO_C \leftarrow \sO_{\PP^3} \leftarrow \oplus_j \sO_{\PP^3}(-j)^{\beta_{1j}} \leftarrow \sG \leftarrow 0
$$ is a resolution by locally free sheaves.
If $C$ is arithmetically Cohen-Macaulay, then $F_3=0$ and $\sG$ splits into a direct sum of line bundles. In this case the ideal $I_C$ is generated by the maximal minors of
$F_1  \leftarrow F_2$ by  the Hilbert-Burch Theorem \ref{HB}. In general we have
 $M \cong \sum_{n \in \ZZ} H^2(\PP^3,\sG(n))$ and $\sum_{n \in \ZZ}H^1(\PP^3,\sG(n))=0$.

The importance of $M$ for liaison comes about as follows. Suppose $f,g \in I_C$ are homogeneous forms  of degree $d$ and $e$ without common factor. Let $X=V(f,g)$ denote the corresponding complete intersection, and let $C'$ be the residual scheme defined by the homogeneous ideal
$I_{C'}=(f,g):I_C$.
The locally free resolutions of $\sO_C$ and  $\sO_{C'}$ are closely related:
Applying $\sExt^2(-,\omega_{\PP^3})$ to the sequence
$$0 \to \sI_{C/X} \to \sO_X \to \sO_C \to 0$$
gives
$$
0 \leftarrow \mathcal Ext^2(\sI_{C/X}, \omega_{\PP^3}) \leftarrow \omega_X \leftarrow \omega_C \leftarrow 0.
$$
From $\omega_X \cong \sO_X(d+e-4)$ we conclude $\mathcal Ext^2(\sI_{C/X},\sO_{\PP^3}(-d-e)) \cong \sO_{C'}$ and hence $\sI_{C'/X} \cong \omega_C(-d-e+4)$.

The mapping cone 
$$
  \xymatrix{
        0&\sO_C \ar[l]             &   \sO \ar[l] &
      \bigoplus_j \sO(-j)^{\beta_{1j}}\ar[l] & \sG \ar[l]&0 \ar[l] \\
     0 &\sO_X \ar[l] \ar[u]   &   \sO \ar[l] \ar[u]_\cong  & 
     \sO(-d) \oplus \sO(-e) \ar[l] \ar[u] & \sO(-d-e) \ar[l] \ar[u]& 0  \ar[l] \\
  }
  $$
dualized with $\sHom(-,\sO(-d-e))$ yields the locally free resolution 
$$
0 \to \bigoplus_j \sO(j-d-e)^{\beta_{1j}} \to \sG^*(-d-e)\oplus \sO(-e) \oplus \sO(-d) \to  \sO \to  \sO_{C'} \to 0
$$
In particular one has

\begin{align*}
M_{C'}&= \sum_{n\in \ZZ} H^1(\PP^3, \sI_{C'} (n)) 
 \cong  \sum_{n \in \ZZ} H^1(\PP^3,\sG^*(n-d-e)) \\
& \cong  \sum_{n\in \ZZ} H^2(\PP^3,\sG(d+e-4-n))^* 
 \cong  \Hom_\FF(M_C,\FF)(4-d-e)
\end{align*}
Thus curves, which are related via an even number of liaison steps, have the same Hartshorne-Rao module up to a twist.
Rao's famous result \cite{Rao} says that the even liaison classes are in bijection with finite length graded $S$-modules up to twist.

Reversing the role of $C$ and $C'$, we conclude that the ideal sheaf of $C$ has a locally free resolution
$$ 0 \leftarrow \mathcal  I_C  \leftarrow \sF \oplus \sL_1 \leftarrow \sL_2   \leftarrow 0$$
where $\sL_1=\oplus_\ell \sO(-c_\ell)$ and $\sL_2=\oplus_k \sO(-d_k)$ are direct sums of line bundles, while $\sF$ is a locally free sheaf without line bundle summands. Note that $\sF$ has no $H^2$-cohomology, and its $H^1$-cohomology 
$$\sum_{n\in \ZZ} H^1(\PP^3,\sF(n)) \cong M$$
is the Hartshorne-Rao module of $C$. The map $\sL_2 \to \sL_1$ coming from a liaison construction might be non-minimal,
in which case one can cancel summands of $\sL_1$ against some summands of $\sL_2$. Below we will work with complexes which arise after such cancellation.

Since there are no line bundle summands and no $H^2$, $\sF$ is  determined by its $H^1$-cohomology:
Consider a minimal presentation of $M$ 
 $$0 \leftarrow M \leftarrow \oplus_i S(-a_i) \leftarrow \oplus _j S(-b_j) \leftarrow N \leftarrow 0 $$
 and the kernel $N$. Then $\sF \cong \widetilde N$ is the associated coherent sheaf of $N$.

 The main difficulty in constructing space curves lies in the construction of $M$. Given $M$ we can find the curve as the cokernel of an homomorphism
 $\varphi \in \Hom(\sL_2,  \sF \oplus \sL_1)$, which we may regard as a vector bundle version of a Hilbert-Burch matrix. Frequently in interesting examples we will have $\sL_1=0$. 
 
 We apply this approach for the case $g=12$ and  $d=13$.  
 
To get an idea about the twists in the resolutions of $M$ and $C$, we use the Hilbert numerators: Suppose the minimal finite free resolution of $M$ is
 $$0 \leftarrow M \lTo F_0 \leftarrow F_1 \leftarrow \ldots \leftarrow F_4 \leftarrow 0$$
 with 
 $$F_i = \oplus S(-j)^{\beta_{ij}}$$
 then the Hilbert series
 $$H_M(t) =\sum_n \dim M_n t^n = \frac{\sum_{i,j}  (-1)^i \beta_{ij} t^j}{(1-t)^4}.$$ 
 Assuming the open condition that $C$ has maximal rank, i.e. 
 $$H^0(\PP^3,\sO(n)) \to H^0(C, L^n)$$ 
 is of maximal rank for all $n$ then
 $H_M(t)=5t^2+8t^3+6t^4$
 and the Hilbert numerator is
 $$H_M(t)(1-t)^4=5t^2-12t^3+4t^4+4t^5+9t^6-16t^{10}+6t^{11}$$
 If $M$ has a natural resolution, which means that for each $j$ at most one $\beta_{ij}$  is non-zero,
 then  $M$ has a Betti table $( \beta_{i,i+j})$
 \begin{verbatim}
            0  1  2  3 4
     total: 5 12 17 16 6
         2: 5 12  4  . .
         3: .  .  4  . .
         4: .  .  9 16 6
\end{verbatim}
 in \textit{Macaulay2} notation. Having natural resolution is another open condition. If  $S_C=S/I_C$ has natural resolution as well, then its Betti table will  be
\begin{verbatim}
            0  1  2 3
     total: 1 11 16 6
         0: 1  .  . .
         1: .  .  . .
         2: .  .  . .
         3: .  .  . .
         4: .  2  . .
         5: .  9 16 6
\end{verbatim}
Comparing these tables, we conclude that $\sL_1=0$, $\rank \sF = 12-5=7$ and $\sL_2= \sO(-4)^4 \oplus \sO(-5)^2$ with $\rank \sL_2 =6$ will hold for an open set of curves. Of course,
we have to prove that the set is non-empty and contains smooth curves.\medskip

\noindent
To construct a desired $M$, we start with the submatrix $\psi$ defining $ S^{12}(-3) {\leftarrow}S^4(-4)$ in $F_1 \leftarrow F_2$ which we  choose randomly.
Since $12{3+1\choose 3}-4{3+2\choose 3}=8 >5$, the kernel $\ker(\psi^t: S^{12}(3) \to S^4(4))$ has at least $8$ linearly independent
elements of degree $-2$, and we choose $5$ random linear combinations of these elements to get the transpose of the presentation matrix $F_0 \leftarrow F_1$ of $M$.
\begin{verbatim}
  i35 : FF=ZZ/10007
  i36 : S=FF[x_0..x_3]
  i37 : psi=random(S^{12:-3},S^{4:-4}) -- the submatrix
  i38 : betti(syzpsit=syz transpose psi)
  i39 : M=coker transpose(syzpsit*random(source syzpsit,S^{5:2}));
  i40 : F= res M -- free resolution of the desired Hartshorne-Rao module
  i41 : betti F
  i42 : L2=S^{4:-4,2:-5}
  i43 : phi= F.dd_2*random(F_2,L2);
  i44 : betti(syzphit=syz transpose phi)
  i45 : IC=ideal mingens ideal( transpose syzphit_{5}*F.dd_2);
  i46 : betti res IC -- free resolution of C   
  i47 : codim IC, degree IC, genus IC
  o47 = (2, 13, 12)
  o47 : Sequence
\end{verbatim}
Next we check smoothness. Since $C\subset \PP^3$ has codimension $2$ and is locally determinantal,
it is unmixed, and we can apply the jacobian criterion.
\begin{verbatim}
  i48 : singC=IC+minors(2,jacobian IC);
  i49 : codim singC==4
  o49 = true
\end{verbatim}
Thus the curve $C$ is smooth. We conclude that  the Hilbert scheme $\operatorname{Hilb}_{d,g}(\PP^3)$ of curves of degree $d=13$ and $g=12$ has a unirational component $H$ defined over $\QQ$, whose general element is smooth. 

As a corollary we get the  unirationality of $\gM_{12}$ from this if we verify that the fiber of the rational map $H \dasharrow \gM_{12}$ has the expected dimension $\dim PGL(4) + \rho(g,r,d)=15+4$ at the given point. Now this
holds if the Petri map 
$H^0(\omega_C(-1)) \otimes H^0( \sO(1)) \to H^0(\omega_C)$ is injective
where $\omega_C= \mathcal Ext^2(\sO_C,\omega_{\PP^3}) \cong \mathcal Ext^1(\sI_C, \sO_{\PP^3}(-4))$ denotes the canonical bundle on $C$.

\begin{verbatim}
  i50 : betti Ext^1(IC,S^{-4})
               0  1
  o50 = total: 6 12
           -1: 2  .
            0: 4 12
\end{verbatim}
shows that there are no linear relation among the two generators of 
$$\Gamma_*(\omega_C)=\sum_{n \in \ZZ} H^0(\PP^3,\omega_C(n))=Ext^1_S(I_C,S(-4))
$$
in degree $-1$.
Thus $H$ dominates and $\gM_{12}$ is unirational. As a corollary of our construction we obtain that the Hurwitz scheme $\mathcal H_{k,g}$ of $k$-gonal curves of genus $g$ is unirational for $(k,g)=(9,12)$.
Indeed $\omega_C(-1)$ is a line bundle of degree $k=22-13=9$.

The case $d=13, g=12$ is actually not used in Sernesi's proof of the unirationality of $\gM_{12}$. Chang, Ran and Sernesi use the cases $(d,g)=(10,11), (12,12)$ and $(13,13)$, all of which can be treated similarly as the case above. I took the case $(d,g)=(13,12)$ because it illustrates the difficulty in constructing $M$ well, and because being different, it implies some minor  new results.
Computer algebra simplifies the cumbersome proof that the approach leads to smooth curves in all of theses cases. For more details and an implementation we refer to the \textit{Macaulay2} package RandomCurves. 

There are 65 values for $(d,g)$ such that there are possibly non-degenerate maximal rank curves with natural resolution
such that the Hartshorne-Rao module has diameter $\le 3$ and natural resolution as well. 
For 60 of these values one can establish the existence of a unirational component in the Hilbert scheme by the methods above.

To prove the unirationality of $\gM_g$ for $g \ge 14$ via space curves, leads to Hartshorne-Rao modules which have diameter $\ge 4$, i.e., modules, which are nonzero in at least 4 different degrees. The construction of such modules is substantially more difficult.

\section{Verra's proof of the unrationality of $\gM_{14} $}\label{Verra}

Verra's idea \cite{Ve} is of beautiful simplicity. Consider a general curve of genus $g=14$ and a pencil $| D |$ of minimal degree which is $\deg D = 8$. The Serre dual linear system $| K-D |$ embeds
$C$ into $\PP^6$ as curve of degree $18$ with expected syzygies 
\begin{verbatim}
             0  1  2  3  4 5
      total: 1 13 45 56 25 2
          0: 1  .  .  .  . .
          1: .  5  .  .  . .
          2: .  8 45 56 25 .
          3: .  .  .  .  . 2
\end{verbatim}
In particular, $C \subset \PP^6$ lies on 5 quadrics. The intersection of the quadrics should have codimension 5 and degree $2^5=32$, thus should equal  $C \cup C'$ where the residual  curve $C'$ has smaller degree $32-18=14$ than $C$ and hence also smaller genus, $g(C')=8$ as it turns out. 

(Indeed the dualizing sheaf of the complete intersection is
$\omega_{C \cup C'} \cong \sO_{C \cup C'}(5\cdot 2-7)\cong \sO_{C\cup C'}(3)$, and we obtain arithmetic genus $p_a(C\cup C')=2^5\cdot 3/2+1=49$. 
Assuming only nodes as intersection points we have 
$\omega_{C \cup C'} \otimes \mathcal  O_C\cong \omega_C (C\cap C')$ and get $\deg (C \cap C')=3\cdot 18-26 = 28$  intersection points. The formula
$p_a(C\cup C')=p_a(C)+p_a(C')+\deg (C\cap C') -1$ finally gives $p_a(C')=49-14-28+1=8$.)
There is no reason to expect anything else than, that $C'$ is a smooth curve.

By a famous result by Mukai \cite{Mu}, every general canonical curve of genus $8$ is obtained as transversal intersection
$$ C' = \mathbb G(2,6) \cap \PP^7 \subset \PP^{14}$$
of the Grassmannian $ \mathbb G(2,6)\subset  \PP^{14}$ in its Pl\"ucker embedding. If we choose
8 general points  on $\mathbb G(2,6)$ and $\PP^7$ as their span, then we get a genus 8 curve together with 8 points. We group these 8 points into two divisors $D_1=p_1+\ldots+p_4$ and 
$D_2=p_5+\ldots+p_8$ of degree 4. Then $\mid K_{C'}+D_1-D_2\mid$ is a general linear system of degree 14 on $C'$, and re-embedding $C' \hookrightarrow \PP^6$ leads to a curve with expected syzygies
\begin{verbatim}
             0 1  2  3  4 5
      total: 1 7 35 56 35 8
          0: 1 .  .  .  . .
          1: . 7  .  .  . .
          2: . . 35 56 35 8
\end{verbatim}

The \textit{Macaulay2} code is now straightforward. First we construct $C'$ in its canonical embedding.
 \begin{verbatim}
  i51 : randomCurveOfGenus8With8Points = R ->(
        --Input: R a polynomial ring in 8 variables, 
        --Output: a pair of an ideal of a canonical curve C
        --        together with a list of ideals of 8 points
        --Method: Mukai's structure theorem on genus 8 curves.
        --  Note that the curves have general Clifford index. 
        FF:=coefficientRing R;
        p:=symbol p;
       -- coordinate ring of the Plucker space:
        P:=FF[flatten apply(6,j->apply(j,i->p_(i,j)))]; 
        skewMatrix:=matrix apply(6,i->apply(6,j->if i<j then 
            p_(i,j) else if i>j then -p_(j,i) else 0_P));
        -- ideal of the Grassmannian G(2,6):
        IGrass:=pfaffians(4,skewMatrix);
        points:=apply(8,k->exteriorPower(2,random(P^2,P^6)));
        ideals:=apply(points,pt->ideal( vars P*(syz pt**P^{-1})));  
        -- linear span of the points:
        L:= ideal (gens intersect ideals)_{0..6};
        phi:=vars P%L; -- coordinates as function on the span
        -- actually the last 8 coordinates represent a basis
        phi2:= matrix{toList(7:0_R)}|vars R; 
        -- matrix for map from R to P/IC
        IC:=ideal (gens IGrass%L); --the ideal of C on the span
        -- obtained as the reduction of the Grassmann equation mod L
        IC2:=ideal mingens substitute(IC,phi2);
        idealsOfPts:=apply(ideals,Ipt->
             ideal mingens ideal sub(gens Ipt%L,phi2));
        (IC2,idealsOfPts))
\end{verbatim}
Building upon Mukai's result, we can construct the desired curve $C'$:
\begin{verbatim}
i52 : randomNormalCurveOfGenus8AndDegree14 = S -> (
      -- Input:  S coordinate ring of P^6
      -- Output: ideal of a curve in P^6               
      x:=symbol x;
      FF:=coefficientRing S
      R:=FF[x_0..x_7];
      (I,points):=randomCurveOfGenus8With8Points(R);
      D1:=intersect apply(4,i->points_i); -- divisors of degree 4 
      D2:=intersect apply(4,i->points_(4+i));
      -- compute the complete linear system |K+D1-D2|, note K=H1       
      H1:=gens D1*random(source gens D1,R^{-1});
      E1:=(I+ideal H1):D1; -- the residual divisor
      L:=mingens ideal(gens intersect(E1,D2)%I); 
      -- the complete linear system
      -- note: all generators of the intersection have degree 2.
      RI:=R/I; -- coordinate ring of C' in P^7       
      phi:=map(RI,S,substitute(L,RI));
      ideal mingens ker phi)
i53 : FF=ZZ/10007;S=FF[x_0..x_6];
i55 : I=randomNormalCurveOfGenus8AndDegree14(S);
i56 : betti res I
\end{verbatim}
Finally, we get curves of genus 14 with
\begin{verbatim}
i57: randomCurveOfGenus14=method(TypicalValue=>Ideal,
          Options=>{Certify=>false})
      -- The default value of the option Certify is false, because 
      -- certifying smoothness is expensive
i58 : randomCurveOfGenus14(PolynomialRing) :=opt ->( S-> (
      -- Input: S PolynomialRing in 7 variables
      -- Output: ideal of a curve of genus 14
      -- Method: Verra's proof of the unirationality of M_14
      IC':=randomNormalCurveOfGenus8AndDegree14(S);
      -- Choose a complete intersection:
      CI:=ideal (gens IC'*random(source gens IC',S^{5:-2}));
      IC:=CI:IC'; -- the desired residual curve
      if not opt.Certified then return IC;
      if not (degree IC ==18 and codim IC == 5 and genus IC ==14) 
           then return nil;
      someMinors :=minors(5, jacobian CI);
      singCI:=CI+someMinors;
      if not (degree singCI==28 and codim singCI==6) 
          then return nulll;
      someMoreMinors:=minors(5, jacobian (gens IC)_{0..3,5});
      singC:=singCI+someMoreMinors;
      if codim singC == 7 then return IC else return nil))
i59 : time betti(J=randomCurveOfGenus14(S))
i60 : time betti(J=randomCurveOfGenus14(S,Certified=>true))
i61 : betti res J
\end{verbatim}
To deduce from these computation the unirationality of $\gM_{14}$, we have to prove again that the  Petri map is injective. By the Betti numbers of the free resolution of $S_C$ we see that there is is no linear relation 
among the two generators of $\omega_C(-1)$. Thus the family is dominant because the conditions that $C$ has expected Betti numbers  and that $C'$ is smooth are open. 

Finally, we remark that, comparing the syzygies of $\sO_C$, $\sO_{C \cup C'}$ and $\sO_{C'}$ via liaison
as outlined for space curves in Section \ref{spaceCurves}, 
we see that
the Petri map of $(C,\sO_C(1))$ is injective, if and only if $I_{C'}$ is generated by quadrics. Indeed, suppose the Betti table of $C'$ is 
\begin{verbatim}
             0 1  2  3  4 5
          0: 1 .  .  .  . .
          1: . 7  x  .  . .
          2: . x 35 56 35 8
\end{verbatim}
where we assume that we need $x$ cubic generators of the ideal of $I_{C'}$.
Then the Betti table of the mapping cone 
with the Koszul complex resolving $\sO_{C \cup C'}$ is
\begin{verbatim}
             0 1  2  3  4 5 6
         -1: . 1  .  .  . . .               
          0: 1 .  5  .  . . .
          1: . 7  x 10  . . .
          2: . x 35 56 45 8 .  
          3: . .  .  .  . 5 .
          4: . .  .  .  . . 1
\end{verbatim}
Minimalizing the dual complex, leads to following Betti table of $C$:
\begin{verbatim}
             0  1  2  3  4 5
          0: 1  .  .  .  . .
          1: .  5  .  .  . .
          2: .  8 45 56 25 x
          3: .  .  .  .  x 2
\end{verbatim}

\section{Minimal Resolution Conjectures and Koszul Divisors}

A graded $S$-module $M$ is said to satisfy the minimal resolution conjecture (MRC) (or is said to have expected syzygies) if the Betti numbers $\beta_{ij}$ of the minimal free resolution
$$0 \lTo M \lTo F_0 \lTo F_1 \lTo \ldots F_c \lTo 0$$
where $F_i= \oplus_j S(-j)^{\beta_{ij}}$ satisfy the following: For each internal degree $j$ at most one of the numbers $\beta_{ij} \not=0$.
In other words, $M$ satisfies MRC if one can nearly  read off the Betti table from the Hilbert numerator
of the Hilbert series $H_M$. 

We say that a module $M$ has a pure resolution, if for  each homological degree $i$ at most one $\beta_{ij} \not=0$.  Betti tables of Cohen-Macaulay modules with pure resolutions play a special role since they span the extremal rays in the Boij-S\"oderberg cone of all Betti tables, see \cite{ES,BS}. Moreover, if $M$
has a pure resolution, then $M$ satisfies the MRC, see \cite{ES}.

We say that the MRC holds generically on a component $H$ of a Hilbert scheme $\operatorname{Hilb}_{p(t)}(\PP^n)$ if it is satisfied for the coordinate ring $S_X$ for  $X \in U \subset H \subset \operatorname{Hilb}_{p(t)}(\PP^n)$ in a dense open subset $U$ of $H$. Note that the Hilbert function is constant on an open set $U'$ of $H$, so this makes sense.
Since Betti numbers are upper semi-continuous in flat families with constant Hilbert function, there is a smallest
possible Betti table for $X \in U'$, and we can ask what this table is. If the MRC is satisfied generically on $H$, then in a sense, the question, what is the generic Betti table, has a trivial answer. 

The MRC has been established in various cases, the most famous one is  the following.

\begin{theorem}[Voisin, 2005 \cite{Vo}, Generic Green's Conjecture \cite{Gr}]  A general canonical curve of genus $g\ge 3$ over a ground field of characteristic 0 satisfies the MRC.
\end{theorem}

In more concrete terms, this conjecture says for example, that a general canonical curve of genus $g=15$ has the  Betti table
\begin{verbatim}
       0  1   2    3    4    5    6    7    8    9   10  11 12 13            
total: 1 78 560 2002 4368 6006 4576 4576 6006 4368 2002 560 78  1
    0: 1  .   .    .    .    .    .    .    .    .    .   .  .  .
    1: . 78 560 2002 4368 6006 4576    .    .    .    .   .  .  .
    2: .  .   .    .    .    .    . 4576 6006 4368 2002 560 78  .
    3: .  .   .    .    .    .    .    .    .    .    .   .  .  1
\end{verbatim}
Note that the symmetry of the table comes from the fact that the homogeneous coordinate ring of a canonical curve is Gorenstein. Before Voisin's result, $g=15$ was for a long time the bound of how far a computer algebra verification of the generic Green Conjecture was feasible, the main problem being the memory requirements in the syzygy computation. 

The full Green conjecture is the following
\begin{conjecture}[Green's Conjecture \cite{Gr}] Let $C$ be a smooth projective curve of genus $g \ge 3$ over a field of characteristic zero. The canonical ring $R_C= \sum_n H^0(C, \omega_C^{\otimes n})$ of $C$
as an $S=Sym H^0(C,\omega_C)$-module has vanishing Betti numbers
$\beta_{i,i+2}=0$ if and only if $i< \Cliff(C)$.
\end{conjecture}
\noindent
Here  the Clifford index is defined as
$$\Cliff(C) = \min \{ \deg L-2(h^0(L)-1) \mid h^0(L), h^1(L) \ge 2 \}.$$
 This conjecture generalizes the
Noether-Petri-Babbage theorem which are the cases $i=0$ and $i=1$ respectively. Its known for $i=2$ in full generality, see \cite{Vo1,Sch2}.
The direction $\Cliff(C) \le i \Rightarrow \beta_{i,i+2} \not=0$ was established by Green and Lazarsfeld in 
the appendix to \cite{Gr}

Green's conjecture is known to be wrong over fields of finite characteristic. The first cases are $g=7$ in characteristic 2 and genus $g=9$ in characteristic 3, see \cite{Sch1,Sch3}. 

Frequently, Green's Conjecture is formulated with the property $N_p$ of Green-Lazarsfeld \cite{GL}:
A subscheme $X \subset \PP^r$ satisfies the property $N_p$ if $\beta_{ij} = 0$ for all $j \ge i+2$ and $i \le p$, in other words if the first $p$ steps in the free resolution are as simple as possible.
An equivalent formulation of Green's Conjecture is that a canonical curve satisfies property $N_p$ iff $p < \Cliff(C)$.

If the coordinate ring $S_X$ for general $X \in U' \subset H \subset \operatorname{Hilb}_{p(t)}(\PP^r)$ is pure, then 
$$
\sZ= \{ X \in U' \mid S_X \hbox{ does not satisfy MRC} \}
$$
is a divisor, a so-called Koszul divisor, because  in principle it can be computed with the Koszul complex and Koszul cohomology groups.

For odd genus $g=2k-1$ Hirschowitz and Ramanan computed the class of the corresponding
divisor in $\gM_g$.

\begin{theorem}[Hirschowitz-Ramanan \cite{HR}]
If the generic Green conjecture holds for odd genus $g=2k-1$, then the Koszul divisor
$$
\sZ=\{ C \in \gM_g \mid R_C \hbox{ does not satisfy MRC } \} = (k-1)\gM_{g,k}^1
$$ 
where $\gM_{g,k}^1= \{ C \in \gM_g \mid \exists \, g^1_k \hbox{ on } C\}$
is the Brill-Noether divisor.
\end{theorem}

The proof is based on a divisor class computation of $\sZ$ and $\gM_{g,k}^1$ on a partial compactification of $\gM_g$ and the fact that
$\sZ -(k-1)\gM_{g,k}^1$ is effective. The coefficient $k-1$ is explained by the fact that a curve with a $g^1_k$ has $\beta_{k-1,k}=\beta_{k-2,k} \ge k-1$ and equality holds for a general $k$-gonal curve of genus $g=2k-1$ as a consequence of the Hirschowitz-Ramanan Theorem.

Based on the Hirschowitz-Ramanan and Voisin's Theorem, Aprodu and Farkas \cite{AF} established that Green's Conjecture holds for smooth curves on arbitrary K3 surfaces. 

One can hope that there are various further interesting Koszul divisors. To get divisors on $\gM_g$ we consider the case that $g=(r+1)s$, $d=r(s+1)$ and ask whether a normal curve of degree $d$, genus $g$ and speciality $h^1(L)=s$ can have a pure resolution. This is a purely numerical condition on
$r$ and $s$.

\begin{conjecture}[Farkas, \cite{Fa1}] Let $p\ge 0$ and $s \ge 2$ be integers. Set $r=(s+1)(p+2)-2, \, g=(r+1)s$ and $d=r(s+1)$. A general smooth normal curve $C \subset \PP^r$ of genus $g$, degree $d$ and speciality $h^1(\sO_C(1))=s$ has a pure resolution, equivalently it satisfies $N_p$. \end{conjecture}

\noindent If the conjecture is true, then the divisor 
$$\sZ_{p,s}=\{ C \in \gM_g \mid \exists g^r_d \hbox{ which  does not satisfy property } N_p \}$$
gives counterexamples to the slope conjecture, see \cite{AF1}. For small cases Farkas verified the conjecture computationally in the spirit of the computations below.

Turning to non-special curves, we have

\begin{conjecture}[Farkas, \cite{FL}] 
A general smooth normal curve $C \subset \PP^r$ of odd genus $g=2p+3$ and degree $d=2g$ has a pure resolution, equivalently satisfies $N_p$.
A general smooth normal curve of even genus $g=2p+6$ and degree $2g-2$ has a pure resolution, equivalently satisfies $N_p$.
\end{conjecture} 

For the rest of this section we return to the task to verify these conjectures computationally in a few cases. The case of genus $g=8$ and degree $d=14$ has been established along the proof of
Verra's Theorem in Section \ref{Verra}. We treat the case of genus $g=7$ and degree $d=14$. This case was settled in \cite{Fa3} using reducible curves. Here we use smooth curves.

Consider for $g=7$ a random plane model of degree $d=7$ with $\delta=8$ nodes.

\begin{verbatim}
  i62: FF=ZZ/10007
  i63: R=FF[x_0..x_2]
  i64: g=7
  i65: delta=binomial(6,2)-7
  i66: J=randomDistinctPlanePoints(delta,R)
  i67: betti res J
  i68: betti (J2=saturate(J^2))
  i69: C=ideal (gens J2*random(source gens J2,R^{-7}))
\end{verbatim}

To find a general divisor of degree $2g=14$ on $C$ we note that over a large finite field geometrically irreducible varieties have always a lot of $\FF$-rational points. We can frequently find points as one of the ideal theoretic components  of the intersection with a random complementary dimensional linear subspace.

\begin{verbatim}
  i70: decompose(C+ideal random(1,R))
  i71: apply(decompose(C+ideal random(1,R)),c->degree c)
  i72: time tally apply(1000,i->
           apply(decompose(C+ideal random(1,R)),c->degree c))	  
     -- used 9.62383 seconds
  o72 = Tally{{1, 1, 1, 1, 1, 2} => 5}
            {1, 1, 1, 1, 3} => 18
            {1, 1, 1, 2, 2} => 21
            {1, 1, 1, 4} => 47
            {1, 1, 2, 3} => 84
            {1, 1, 5} => 91
            {1, 2, 2, 2} => 19
            {1, 2, 4} => 120
            {1, 3, 3} => 58
            {1, 6} => 167
            {2, 2, 3} => 40
            {2, 5} => 94
            {3, 4} => 77
            {7} => 159
  o72 : Tally

\end{verbatim}
Thus the following function will provide points.
\begin{verbatim}
  i73: randomFFRationalPoint=method()
  i74: randomFFRationalPoint(Ideal):=I->(
       --Input: I ideal of a projective variety X
       --Output: ideal of a FF-rational point of X
       R:=ring I;
      if char R == 0 then error "expected a finite ground field";
      if not class R === PolynomialRing then 
          error "expected an ideal in a polynomial ring";
      n:=dim I-1;
      if n==0 then error "expected a positive dimensional scheme";
      trial:=1;
      while (
          H:=ideal random(R^1,R^{n:-1});
          pts:=decompose (I+H);
         pts1:=select(pts,pt-> degree pt==1 and dim pt ==1);
        #pts1<1 ) do (trial=trial+1);
      pts1_0)
 i75: randomFFRationalPoint(C)
\end{verbatim}

This allows to get a random effective divisor of degree $14$ on $C$ with all points $\FF$-rational,
for which we verify $N_2$:

\begin{verbatim}
  i76: points=apply(2*g,i->randomFFRationalPoint(C))
  i77: D=intersect points -- effective divisor of degree 14
  i78: degree D
  o78 = 14
  i79: DJ=intersect(D,J)
  i80: degree DJ==degree D + degree J
  o80 = true
  i81: betti DJ
  i82: H=ideal(gens DJ*random( source gens DJ, R^{-6}))+C
  i83: E=(H:J2):D -- the residual divisor
  i84: degree E +degree D + 2*degree J == 6*7
  o84 = true
  i85: L=mingens ideal ((gens truncate(6,intersect(E,J))) %C)
         -- the complete linear series |D|
  i86: RC=R/C
  i87: S=FF[y_0..y_7]
  i88: phi=map(RC,S,sub(L,RC))
  i89: I=ideal mingens ker phi;-- C re-embedded
  i90: (dim I,degree I, genus I) == (2,14,7)
  i91: time betti res I
               0  1  2  3  4  5 6
  o91 = total: 1 14 28 56 70 36 7
            0: 1  .  .  .  .  . .
            1: . 14 28  .  .  . .
            2: .  .  . 56 70 36 7
\end{verbatim}

The computation above is the reduction mod $p$ of the computation for a curve together with a divisor
over some open part of $\Spec O_K$ of a number field $K$. We can bound the degree $[K:\QQ] \le 7^{14}$ of the number field, a bound on its discriminant seems out of reach. 

Hence as before we can conclude that Farkas' syzygy conjecture holds for curves of genus 7 and degree 14 over  fields of characteristic zero.
Using Mukai's description of curves of genus 7 \cite{Mu}, we could find an example defined over $\mathbb Q$. I chose to present the example above because I wanted to illustrate the trick to get points
in $\gM_{g,n}$ from points in $\gM_g$ over finite fields. 

It is not difficult to  use similar constructions to verify Farkas' conjecture for other small $g$.

It is known that the MRC is not satisfied in a number of important cases, notably for non-special normal curves of genus $g\ge4$ and large degree \cite{Gr}.

\vbox{\noindent Author Address:\par
\smallskip

\noindent{Frank-Olaf Schreyer}\par
\noindent{Mathematik und Informatik, Universit\"at des Saarlandes, Campus E2 4, 
D-66123 Saarbr\"ucken, Germany}\par
\noindent{schreyer@math.uni-sb.de}\par
}

\end{document}